\documentclass[reqno,11pt]{amsart}

\usepackage{xcolor}
\usepackage{graphicx}
\usepackage[hidelinks]{hyperref}
\usepackage{amscd}

\usepackage{amsmath}
\usepackage{amsthm}
\usepackage{amssymb}
\usepackage{latexsym,array}
\usepackage{amsfonts}
\usepackage{shadow}
\usepackage{amsbsy}
\usepackage{amssymb}
\usepackage{dsfont}
\usepackage{mathtools}
\usepackage{enumitem}
\usepackage{doi}

\usepackage[notref, notcite, final]{showkeys}
\usepackage{color}
\usepackage{graphicx}
\usepackage[hidelinks]{hyperref}
\usepackage{cleveref}
\usepackage{fullpage}
\usepackage{comment}
\usepackage{tikz-cd}
\usepackage{tikz}
\usepackage{float}

\usepackage{dsfont}

\usepackage{cleveref}

\numberwithin{equation}{section}

\newtheorem{theorem}{Theorem}[section]

\newtheorem{corollary}[theorem]{Corollary}
\newtheorem{proposition}[theorem]{Proposition}

\newtheorem{exs}[theorem]{Examples}

\theoremstyle{definition}
\newtheorem{remark}[theorem]{{\bf Remark}}
\newtheorem{definition}[theorem]{Definition}

\title[Gaussian RBF-kernels via Fock spaces: quaternionic and several complex variables settings]{Gaussian RBF-kernels via Fock spaces: quaternionic and several complex variables settings}
\date{}

\author[Antonino De Martino]{Antonino De Martino}
\address{(ADM)
	Schmid College of Science and Technology,  Chapman University One University Drive Orange\\
	California 92866,
	USA
}\email{antonino.demartino@polimi.it}

\author[Kamal Diki]{Kamal Diki}
\address{(KD)
	Schmid College of Science and Technology,  Chapman University One University Drive Orange\\
	California 92866,
	USA
}\email{diki@chapman.edu}

\begin{document}
\maketitle
	\begin{abstract}
In this paper we study two extensions of the complex-valued Gaussian radial basis function (RBF) kernel and discuss their connections with Fock spaces in two different settings. First, we introduce the quaternonic Gaussian RBF kernel constructed using the theory of slice hyperholomorphic functions. Then, we consider the case of Gaussian RBF kernels in several complex variables.
\end{abstract}

\noindent AMS Classification: 30G35, 30H20, 44A15, 46E20

\noindent Keywords: Reproducing kernel Hilbert spaces, Gaussian RBF kernel, RBF spaces quaternions, Fock space, several complex variables, Segal-Bargmann transform 
\section{Introduction and preliminary results}
Kernels and reproducing kernel Hilbert spaces (RKHSs) appear in different areas of mathematics such as complex analysis, operator theory and Schur analysis. These two concepts are fundamental in machine learning (ML), particularly in the context of kernel methods, such as Support Vector Machines (SVMs) and Principal Component Analysis (PCA). In fact, kernels are used to measure the similarity between data points, allowing the construction of nonlinear decision boundaries in SVMs and other kernel-based algorithms. For further discussions on kernels, reproducing kernel Hilbert spaces, and their various applications in machine learning we refer the reader to \cite{SHS, SC2008, vert2004primer} and the references therein. Some of the most popular kernels used in SVMs are the so-called Gaussian RBF kernels which we plan to study further in this paper. In addition to that, it is important to note that kernels are used to define coherent states in quantum mechanics, see \cite{Gaz}. In this section we start by reviewing some basic notions on kernels and RKHSs:
\begin{definition}
 Let $\mathcal{X}$ be a non-empty set and $\mathbb{K}=\mathbb{R}$ or $\mathbb{C}$.
\begin{itemize}
\item[i)] A function $k:\mathcal{X}\times \mathcal{X}\longrightarrow \mathbb{K}$ is called a kernel on $\mathcal{X}$ if there exists a $\mathbb{K}$-Hilbert space $\mathcal{H}$ and a map $\Phi:\mathcal{X}\longrightarrow \mathcal{H}$  such that for every $x,y\in \mathcal{X}$ we have
$$k(x,y)=\langle \Phi(y), \Phi(x) \rangle_{\mathcal{H}}.$$
\item[ii)] The Hilbert space $\mathcal{H}$ is called a feature space 
and the map $\Phi$ is called a feature map for the kernel $k$.
\end{itemize}
\end{definition}

The feature map $\Phi$ and the feature space $\mathcal{H}$ are not unique, only the kernel $k$ is unique. On the Euclidean space $\mathbb{R}^d$ we can consider the Euclidean inner product and its associated metric, given by 
$$
 \langle x,y \rangle_{2,d}=\displaystyle \sum_{\ell=1}^{d}x_\ell y_\ell, \quad\displaystyle d(x,y)=||x-y||_{2,d}=\sqrt{\sum_{\ell=1}^{d}|x_\ell-y_\ell|^2},
$$
for every $x=(x_1,...,x_d), y=(y_1,...,y_d)\in\mathbb{R}^d$. Based on the above notations we give some examples of well-known real-valued kernels:

\begin{exs}
For $x=(x_1,...,x_d);y=(y_1,...,y_d)\in\mathbb{R}^d$ we present the following kernels: 

\begin{enumerate}
\item Polynomial kernel: $k(x,y)=(1+\langle x,y \rangle_{2,d})^m, \quad m\geq 1.$
\item Exponential kernel: $k(x,y)=\exp(\langle x,y \rangle_{2,d})$.
\item Radial basis function kernels: $k(x,y)=\varphi(||x-y||)$.
\item Gaussian RBF kernels: we will discuss this example in details in the next sections.
\end{enumerate}
\end{exs}

To each kernel function $k$ we associate a Hilbert space which is generally called \textit{"reproducing kernel Hilbert space"}. This space can be defined as follows:
\begin{definition}
A Hilbert space $\mathcal{H}$ of functions defined on an open set $\Omega$ is called a reproducing kernel Hilbert space (RKHS) if the point evaluations
$$\Lambda_w: f\longmapsto f(w),\quad  w\in\Omega$$ are bounded.
\end{definition}
\begin{remark}
Let $\mathcal{H}$ be a RKHS. Then, by Riesz representation theorem there exists a uniquely determined function $K(z, w)$ defined on $\Omega \times \Omega$ satisfying:

\begin{itemize}
\item[i)] For every $w\in\Omega$, the function $$K_w: z\longmapsto K(z,w)$$ belongs to $\mathcal{H}$.
\item[ii)] Reproducing kernel property: for every $f\in\mathcal{H}$ and $w\in\Omega$, we have
$$\langle f, K_w\rangle_\mathcal{H}=f(w).$$
\end{itemize}
The function $K(z,w)$ is positive definite and is called the reproducing kernel of $\mathcal{H}$.
\end{remark}
  Conversely, we have the following fundamental result 
\begin{theorem}[Moore-Aronszajn theorem]
Associated to a function $K(z,w)$ positive definite on a set $\Omega$ is uniquely determined a Hilbert space $\mathcal{H}(K)$, whose elements are functions on $\Omega$, and with reproducing kernel $K(z, w)$.
\end{theorem}
\begin{remark}
The main take-away message from the previous discussion is the following:
 \begin{center}
     Positive definite functions $\longleftrightarrow$ Kernels $\longleftrightarrow$  Reproducing kernel Hilbert spaces.
    \end{center}
\end{remark}

In this paper we focus on extending the complex Gaussian RBF kernel and relate it to the Fock spaces in the quaternionic and several complex variables settings. To this end, let us review the RBF kernel in one complex variable which was introduced in \cite{SDC} (see also \cite{SC2008}) as follows
\begin{definition}
 Let $\gamma >0$, the complex-valued Gaussian RBF kernel (in one dimension) is denoted by $K_\gamma:\mathbb{C}\times \mathbb{C}\longrightarrow \mathbb{C}$ and can be expressed as follows
$$K_\gamma(z,w)=\exp\left(-\dfrac{(z-\overline{w})^2}{\gamma^2}\right), \quad \forall z,w\in\mathbb{C}. $$
\end{definition}
The associated complex Gaussian RBF space was also introduced in \cite{SDC} (see also \cite{SC2008}) and it is given by  \begin{definition}
Let $\gamma>0$, an entire function $f:\mathbb{C}\longrightarrow \mathbb{C}$ belongs to the RBF space, denoted by $\mathcal{H}_{\gamma}^{RBF}(\mathbb{C})$ (or simply $\mathcal{H}_\gamma)$ if we have
$$
\displaystyle ||f||_{\mathcal{H}_\gamma}^2:=\left(\dfrac{2}{\pi\gamma^2}\right)\int_{\mathbb{C}}|f(z)|^2\exp\left(\frac{(z-\overline{z})^2}{\gamma^2}\right) dA(z)<\infty,$$
where $dA(z)=dxdy$ is the Lebesgue measure with respect to the variable $z=x+iy$.
\end{definition}

The authors of \cite{ADCS} developed a new approach allowing to study the complex Gaussian RBF kernel via the theory of Fock spaces. In this paper we present two extensions of this approach. In Section 2 we introduce a quaternionic Gaussian RBF kernel and discuss its connection with the slice hyperholomorphic Fock space and quaternionic Segal-Bargmann transform. Then, in Section 3 we treat the case of the Gaussian RBF kernel in several complex variables and use the theory of Fock spaces to present a reproducing kernel property for the RBF kernel. Finally, in Section 4 we summarize and present a scheme showing different possible extensions of the Gaussian RBF kernel to various settings.

\section{Gaussian RBF-kernels via Fock spaces: quaternionic slice hypercomplex case}
In this section we extend some results obtained in \cite{ADCS} to the case of slice hyperholomorphic functions on quaternions. \\ \\
In hypercomplex analysis two theories of functions are prominent: the slice hyperholomorphic and the monogenic one. The main difference between these two theories is that polynomials and power series with quaternionic coefficients to the right (or to the left) are slice hyperholomorphic but they are not monogenic (or Fueter regular in the quaternionic case). The theory of monogenic functions is defined by means of an extension of the Cauchy-Riemann operator in $ \mathbb{R}^4$, see \cite{red}. The match between the two function theories is given by the so-called Fueter theorem, see \cite{CSS3, F}.

\medskip

Nowadays the theory of slice hyperholomorphic functions has several applications: the quaternionic spectral theory on the $S$-spectrum (\cite{6CG,CGK, CSS1}), characteristic operator functions and applications to linear system theory \cite{6COFBook}, Schur analysis (\cite{ACSbook}), new classes of fractional diffusion problems based on fractional powers of quaternionic linear operators(\cite{frac5,6CG,frac3}) and peculiar integral transforms (\cite{DMD, DMD2}).
\\ In this paper we will use the theory of slice hyperholomorphic functions to get our results. In the following subsection we review the main notions of this function theory.
\subsection{Preliminaries}

To make the paper self-contained we recall the algebra of quaternions, that is defined as follows
$$ \mathbb{H}:= \{q=q_0+iq_1+jq_2+kq_3 \, | \quad q_0, q_1,q_2,q_3 \in \mathbb{R}\},$$
where the imaginary units satisfy the following relations
$$ i^2=j^2=k^2=-1 \qquad ij=-ji=k, \, jk=-kj=i, \, ki=-ik=j.$$
The quaternionic conjugate is defined as $ \bar{q}=q_0-iq_1-jq_2-kq_3 $, and it satisfies the following property 
$$ \overline{pq}=\bar{q} \bar{p},$$ 
for $p,q \in \mathbb{H}$. The modulus of a quaternion is defined as $|q|= \sqrt{q_0^2+q_1^2+q_2^2+q_3^2}$. The unit sphere of imaginary units in $ \mathbb{H}$ is defined as
$$ \mathbb{S}:=\{q \in \mathbb{H} \, | \, q^2=-1\}.$$ 
We observe that for some real numbers $x$, $y>0$ and imaginary unit $I \in \mathbb{S}$, any quaternion $q \in \mathbb{H}\setminus \mathbb{R}$ can be written as $q=x+Iy$. For any $I \in \mathbb{S}$ we define $\mathbb{C}_I=\mathbb{R}+I \mathbb{R}$. This is isomorphic to the complex plane $ \mathbb{C}$. We observe that their union give the whole quaternions
$$ \mathbb{H}=\bigcup_{I \in \mathbb{S}} \mathbb{C}_I.$$
In this section we use the following type of functions.

\begin{definition}
A function $f:\mathbb{H} \to \mathbb{H}$ of the form
	$$ f(q)=f(x+Iy)= \alpha(x,y)+I \beta(x,y) \qquad \left(\hbox{resp.} \quad f(q)=f(x+Iy)= \alpha(x,y)+ \beta(x,y)I\right),$$
	is left (resp. right) slice hyperholomorphic if $\alpha$ and $ \beta$ are quaternionic-valued functions and satisfy the so-called "even-odd" conditions i.e.
$$
		\alpha(x,y)=\alpha(x,-y), \qquad \beta(x,y)= - \beta(x,y) \qquad \hbox{for all} \qquad (x,y) \in \mathbb{R}^2.
$$
	Moreover, the functions $ \alpha$ and $ \beta$ have to satisfy the Cauchy-Riemann system
	$$ \partial_{x} \alpha(x,y)- \partial_y \beta(x,y)=0, \quad \hbox{and} \quad \partial_{y} \alpha(x,y)+ \partial_x \beta(x,y)=0.$$
\end{definition} 

The set of left (resp. right) slice hyperholomorphic functions on $\mathbb{H}$ is denoted by $\mathcal{SH}_L(\mathbb{H})$ (resp. $\mathcal{SH}_R(\mathbb{H}))$. If the functions $ \alpha$ and $ \beta$ are real-valued functions, then we are dealing with the subset of intrinsic slice hyperholomorphic functions, denoted by $ \mathcal{N}(\mathbb{H})$.
\\ In \cite{CSS1, CSS2, GS2007} it has been proved that slice hyperholomorphic functions have an expansion in series.

\begin{theorem}
	Let $f$ be a (left) slice hyperholomorphic function. Then for any real point $p_0$ in $\mathbb{H}$, the function $f$ can be represented by power series
	$$f(q)=\sum_{m=0}^{+\infty} (q-p_0)^ma_m$$
	on the ball $B(p_0,R)=\{q\in \mathbb{H}; |q-p_0|<R\}$.
\end{theorem}
The pointwise product of two different slice hyperholomorphic functions is not slice hyperholomorphic. However the product of an intrinsic slice hyperholomorphic function and a slice hyperholomorphic function preserves the slice structure and the slice hyperholomorphicity, see \cite[Thm. 2.1.4]{CGK}.

\begin{theorem}
\label{prod}
If $f \in \mathcal{N}(\mathbb{H})$ and $g \in \mathcal{SH}_L(\mathbb{H})$ then $fg \in  \mathcal{SH}_L(\mathbb{H})$. Similarly, if $f \in \mathcal{SH}_R(\mathbb{H})$ and $g \in \mathcal{N}(\mathbb{H})$ then $fg \in  \mathcal{SH}_R(\mathbb{H})$.
\end{theorem}

\subsection{Slice hyperholomorphic Fock space}

Now, we recall some basic notion of the slice hyperholomorphic Fock space, firstly introduced in \cite{ACSS} and further investigated in \cite{DG}.
\begin{definition}
	Let $\nu>0$ be a real parameter. For a given $I \in \mathbb{S}$ we define the slice hyperholomorphic Fock space as
	$$ \mathcal{F}_{Slice}^{\nu}(\mathbb{H}):= \left \{ f \in \mathcal{SH}_L(\mathbb{H}) \,: \, \left(\frac{\nu}{\pi}\right) \int_{\mathbb{C}_I} |f_I(q)| e^{-\nu|q|^2} d \lambda_I(q)<\infty \right\},$$
\end{definition}
where $f_I=f_{|\mathbb{C}_I}$ and $ d\lambda_I(q)=dxdy$ is the Lebesgue measure with respect to the variable $q=x+Iy$.
The right quaternionic space $ \mathcal{F}_{Slice}^{\nu}(\mathbb{H})$ is endowed with the inner product
$$ \langle f,g \rangle_{\mathcal{F}_{Slice}^{\nu}(\mathbb{H})}=  \left(\frac{\nu}{\pi}\right)  \int_{\mathbb{C}_I} \overline{g_I(q)} f_I(q) e^{- \nu |q|^2} d \lambda_I(q),$$
for $f$, $g \in \mathcal{F}_{Slice}^{\nu}(\mathbb{H})$. The associated norm is defined as
$$ \| f \|^2_{\mathcal{F}_{Slice}^\nu(\mathbb{H})}= \left(\frac{\nu}{\pi}\right)\int_{\mathbb{C}_I} |f_I(q)|^2 e^{- \nu |q|^2} d \lambda_I(q).$$
The space $ \mathcal{F}_{Slice}^\nu(\mathbb{H})$ does not depend on the choice of $I \in \mathbb{S}$. In \cite{ACSS} the authors showed that  the monomials $q^n$ form an orthogonal basis of $ \mathcal{F}_{Slice}^\nu(\mathbb{H})$ with
\begin{equation}
	\label{ort}
	\langle q^m, q^n \rangle_{ \mathcal{F}_{Slice}^\nu(\mathbb{H})}= \frac{m!}{\nu^{m}} \delta_{n,m}.
\end{equation}
The space $ \mathcal{F}_{Slice}^\nu(\mathbb{H})$ is a (right) quaternionic Hilbert space with reproducing kernel given by
$$ K(p,q)=e_{*}^\nu(q\bar{p})= \sum_{n=0}^\infty \frac{\nu^n q^n \bar{p}^n}{n!},$$
see \cite[Thm 3.10]{ACSS}.
The reproducing kernel property can be expressed in terms of the following integral representation
\begin{equation}
	\label{zero}
	f(p)= \langle f, K(p,.) \rangle_{\mathcal{F}_{Slice}^{\nu}(\mathbb{H})}= \left(\frac{\nu}{\pi}\right) \int_{\mathbb{C}_I} \overline{e_{*}^{\nu}(q \bar{p})} f_I(q) e^{- \nu |q|^2} d \lambda_I(q).
\end{equation}
The slice hyperholomorphic Fock space satisfies the following sequential characterization
\begin{equation}
\label{seq}
\mathcal{F}_{Slice}^{\nu}(\mathbb{H})= \left \{f(p)= \sum_{n=0}^\infty p^na_n, \quad \{a_n\}_{n \geq 0} \subset \mathbb{H} \quad | \quad \sum_{n=0}^\infty \frac{k!}{\nu^k}|a_n|^2< \infty \right\},
\end{equation}
see \cite[Prop. 3.11]{ACSS}.
\subsection{Slice hyperholomorphic RBF space}
Inspired from the complex case, see \cite{SDC, SC2008}, we give the following.
\begin{definition}
Let $\gamma>0$, the slice hyperholomorphic RBF space is defined as
$$
\mathcal{H}_{\gamma}^{I}(\mathbb{H}):= \left\{f \in \mathcal{SH}_L(\mathbb{H})\, : \,  \left(\frac{2}{ \pi \gamma^2}\right) \int_{\mathbb{C}_I} |f_I(q)| e^{\frac{(q- \bar{q})^2}{\gamma^2}} d\lambda_I(q) < \infty \right\}, \quad \forall I \in \mathbb{S}.
$$
\end{definition}
The right $ \mathbb{H}$-vector space $ \mathcal{H}_{\gamma}^I(\mathbb{H})$ is endowed with the inner product
$$ \langle f,g \rangle_{\mathcal{H}_{\gamma}^{I}(\mathbb{H})}=\left(\frac{2}{ \pi \gamma^2}\right) \int_{\mathbb{C}_I} \overline{g_I(q)} f_I(q)e^{\frac{(q- \bar{q})^2}{\gamma^2}} d \lambda_I(q) ,$$
for $f$, $g \in \mathcal{H}_{\gamma}^I(\mathbb{H})$. The associated norm is given by
$$ \| f\|_{\mathcal{H}_{\gamma}^{I}(\mathbb{H})}^2 = \left(\frac{2}{ \pi \gamma^2}\right)\int_{\mathbb{C}_I} |f_I(q)|^2e^{\frac{(q- \bar{q})^2}{\gamma^2}} d \lambda_I(q).$$
By following similar arguments performed in \cite[Thm. 3.1]{DG} we can show that the slice hyperholomorphic RBF-kernel does not depend on the choice of the imaginary unit $I \in \mathbb{S}$. Therefore, from now on we will denote the space $\mathcal{H}^I_{\gamma}(\mathbb{H})$ simply by $\mathcal{H}_{\gamma,S}(\mathbb{H})$.

\begin{theorem}
An orthonormal basis of the slice hyperholomorphic RBF space is given by 
\begin{equation}
\label{basis}
e_{n}^{\gamma}(q)= \sqrt{\frac{2^n}{\gamma^{2n} n!}} q^n e^{-\frac{q^2}{\gamma^2}}, \qquad \gamma >0.
\end{equation}
\end{theorem}
\begin{proof}
First of all we compute the following integral
\begin{equation}
\label{inner}
\langle e_n , e_m \rangle_{\mathcal{H}_{\gamma,S}}= \left(\frac{2}{ \pi \gamma^2}\right)\int_{\mathbb{C}_I} \overline{e_m(q)} e_n(q) e^{\frac{(q-\bar{q})^2}{\gamma^2}} d \lambda_I(q).
\end{equation}
By formula \eqref{ort} we get
\begin{eqnarray*}
\langle e_n , e_m \rangle_{\mathcal{H}_{\gamma,S}}&=& \frac{2^{n+1}}{\gamma^{2n+2} \pi  n!} \int_{\mathbb{C}_I} \bar{q}^m q^n e^{- \frac{\bar{q}^2}{\gamma^2}} e^{-\frac{q^2}{\gamma^2}}e^{\frac{(q- \bar{q})^2}{\gamma^2}} d \lambda_I(q)\\
&=& \frac{2^{n}}{\gamma^{2n} n!} \langle q^m, q^n \rangle_{\mathcal{F}_{Slice}^{\frac{2}{\gamma^2}}}\\
&=&  \delta_{n,m}.
\end{eqnarray*}
This implies that $e_n^{\gamma}$ belongs to the slice hyperholomorphic RBF space and that two of them are orthogonal.
\\ Now, we show that $e_n^\gamma(q)$ forms a basis. Let $f \in \mathcal{H}_{\gamma,S}(\mathbb{H})$, then the function $q \mapsto e^{\frac{q^2}{\gamma^2}}f(q)$ is slice hyperholomorphic, since it is a product of a slice hyperholomorphic function and an intrinsic function, see Theorem \ref{prod}. This implies that for $\{a_n\}_{n \geq 0} \subset \mathbb{H}$ we have the following expansion in series
\begin{equation}
\label{f1}
f(q)=\sum_{n=0}^\infty q^n e^{- \frac{q^2}{\gamma^2}}a_n= \sum_{n=0}^\infty \sqrt{\frac{\gamma^{2n} n!}{2^n}} e_n^{\gamma}(q) a_n.
\end{equation}
This means that $e_{n}^{\gamma}(q)$ are generators. Now, we show that $e_{n}^{\gamma}(q)$ are independent. By using formula \eqref{f1} and the expression \eqref{inner} we get
\begin{eqnarray*}
\langle f, e_{n}^{\gamma}(p) \rangle_{\mathcal{H}_{\gamma,S}(\mathbb{H})}&=& \left(\frac{2}{ \pi \gamma^2}\right) \int_{\mathbb{C}_I} \overline{e_n^{\gamma}(q)}f_I(q)e^{\frac{(q- \bar{q})^2}{\gamma^2}} d \lambda_I(q)\\
 &=&\left(\frac{2}{ \pi \gamma^2}\right)\int_{\mathbb{C}_I} \overline{e_n^{\gamma}(q)} \left( \sum_{m=0}^\infty \sqrt{\frac{\gamma^{2m} m!}{2^m}} e_m^{\gamma}(q) a_m\right) e^{\frac{(q- \bar{q})^2}{\gamma^2}} d \lambda_I(q)\\
  &=&\sum_{m=0}^\infty  \sqrt{\frac{\gamma^{2m} m!}{2^m}} \left[\left(\frac{2}{ \pi \gamma^2}\right)\int_{\mathbb{C}_I} \overline{e_n^{\gamma}(q)}  e_m^{\gamma}(q) e^{\frac{(q- \bar{q})^2}{\gamma^2}} \lambda_I(q) \right] a_m \\
&=&\sqrt{\frac{\gamma^{2n} n!}{2^n }} a_n.
\end{eqnarray*}
We can exchange the series with the integral because the series expansion converges uniformly on $|q| \leq  r$. Therefore $\langle f, e_{n}^{\gamma}(p) \rangle_{\mathcal{H}_{\gamma,S}(\mathbb{H})}=0$ if and only if for all $n \in \mathbb{N} $ we have $a_n=0$, i.e. $f=0$.

\end{proof}


Now, we show a relation between the spaces $\mathcal{H}_{\gamma,S}(\mathbb{H})$ and $ \mathcal{F}_{Slice}^{\nu}(\mathbb{H})$ where $\nu=\frac{2}{\gamma^2}$.
\begin{theorem}
\label{one}
Let $\gamma>0$, a left slice hyperholomorphic function $f: \mathbb{H} \to \mathbb{H}$ belongs to the slice hyperholomorphic RBF space $\mathcal{H}_{\gamma,S}(\mathbb{H})$ if and only if there exists a unique function $g$ in the slice hyperholomorphic Fock space $ \mathcal{F}_{Slice}^{\frac{2}{\gamma^2}}(\mathbb{H})$ such that
$$ f(q)=e^{- \frac{q^2}{\gamma^2}} g(q), \qquad \forall q \in \mathbb{H}.$$
Furthermore, an isometric isomorphism between the spaces $\mathcal{H}_{\gamma,S}(\mathbb{H})$ and $ \mathcal{F}_{Slice}^{\frac{2}{\gamma^2}}(\mathbb{H})$ is given by
$$ \mathcal{M}^{\gamma^2}[f](q)= e^{\frac{q^2}{\gamma^2}} f(q), \quad \forall f \in \mathcal{H}_{\gamma,S}(\mathbb{H}), \quad q \in \mathbb{H}.$$
\end{theorem}
\begin{proof}
By hypothesis we know that $g(q)=e^{\frac{q^2}{\gamma^2}} f(q)$. This function is slice hyperholomorphic for every $q \in \mathbb{H}$, since it is a product between a slice hyperholomorphic function and an intrinsic slice hyperholomorphic function $e^{\frac{q^2}{\gamma^2}}$, see Theorem \ref{prod}. Now, we show that $g \in \mathcal{F}_{Slice}^{\frac{2}{\gamma^2}}(\mathbb{H})$. Since $f \in \mathcal{H}_{\gamma,S}(\mathbb{H})$ we have
\begin{eqnarray}
\nonumber
\| g \|_{\mathcal{F}_{Slice}^{\frac{2}{\gamma^2}}(\mathbb{H})}^2&=&  \left(\frac{2}{\gamma^2 \pi}\right)\int_{\mathbb{C}_I} |g_I(q)|^2 e^{- \frac{2}{\gamma^2}|q|^2} d \lambda_I(q)\\
\nonumber
&=& \left(\frac{2}{\gamma^2 \pi}\right) \int_{\mathbb{C}_I} |f_I(q)|^2 e^{\frac{\bar{q}^2}{\gamma^2}+ \frac{q^2}{\gamma^2}- \frac{2}{\gamma^2}|q|^2} d \lambda_I(q)\\
\nonumber
&=&\left(\frac{2}{\gamma^2 \pi}\right) \int_{\mathbb{C}_I} |f_I(q)|^2 e^{\frac{(q- \bar{q})^2}{\gamma^2}} d \lambda_I(q)\\
&=& \| f \|_{\mathcal{H}_{\gamma,S}(\mathbb{H})}^2< \infty.
\end{eqnarray}
The previous computations imply that
$$ \| \mathcal{M}^{\gamma^2}[f] \|_{\mathcal{F}_{Slice}^{\frac{2}{\gamma^2}}(\mathbb{H})}= \|f \|_{\mathcal{H}_{\gamma,S}(\mathbb{H})} \qquad \forall f \in \mathcal{H}_{\gamma,S}(\mathbb{H}).$$
This is an isometry property for the operator $ \mathcal{M}^{\gamma^2}$. Moreover, we have
$$ \mathcal{M}^{\gamma^2} \left(e_{n}^{\gamma}(q)\right)= \sqrt{\frac{2^n}{\gamma^{2n} n!}} q^n, \qquad \forall q \in \mathbb{H}.$$
This means that the operator $ \mathcal{M}^{\gamma^2}$ maps orthonormal basis of the space $ \mathcal{H}_{\gamma,S}(\mathbb{H})$ onto an orthonormal basis of the space  $ \mathcal{F}_{Slice}^{\frac{2}{\gamma^2}}(\mathbb{H})$. Hence the operator $ \mathcal{M}^{\gamma^2}$ is a surjective operator.
\end{proof}
As a corollary of the previous theorem we have the following result.

\begin{corollary}
Let $\gamma>0$, the inverse operator of $\mathcal{M}^{\gamma^2}$ is its adjoint and it is given by
$$ \left( \mathcal{M}^{\gamma^2}\right)^{-1}: \mathcal{F}_{Slice}^{\frac{2}{\gamma^2}}(\mathbb{H}) \to \mathcal{H}_{\gamma,S}(\mathbb{H}).$$
This can be computed in the following way
$$ \left( \mathcal{M}^{\gamma^2}\right)^{-1}=\left(\mathcal{M}^{\gamma^2}\right)^{*}=\mathcal{M}^{-\gamma^2}.$$
\end{corollary}

Now, we compute an explicit expression of the reproducing kernel of the space $ \mathcal{H}_{\gamma,S}(\mathbb{H})$.

\begin{theorem}
\label{two} Let $\gamma>0$. The slice hyperholomorphic RBF space $ \mathcal{H}_{\gamma,S}(\mathbb{H})$ is a reproducing kernel Hilbert space whose kernel is given by
$$ K_{\gamma,S}(q,p)=  e^{- \frac{q^2}{\gamma^2}} e_{*}^{\frac{2}{\gamma^2}}(q \bar{p})e^{- \frac{\bar{p}^2}{\gamma^2}}.$$
The reproducing kernel property is given by the following integral representation
$$
f(p)=  \left(\frac{2}{ \pi \gamma^2}\right) \int_{\mathbb{C}_I} \overline{K_{\gamma,S}(q,p)} f_I(q) e^{\frac{(q- \bar{q})^2}{\gamma^2}} d \lambda_I(q), \quad q \in \mathbb{H}, \qquad \forall f \in \mathcal{H}_{\gamma,S}(\mathbb{H}).
$$
\end{theorem}
\begin{proof}
We start observing that since the function $e^{\frac{(q- \bar{q})^2}{\gamma^2}}$ is real-valued, then it can commute with $f_I(q)$, so we have
\begin{eqnarray*}
 \left(\frac{2}{ \pi \gamma^2}\right) \int_{\mathbb{C}_I} \overline{K_{\gamma,S}(q,p)} e^{\frac{(q- \bar{q})^2}{\gamma^2}} f_I(q) d \lambda_I(q)&=&\left(\frac{2}{\pi \gamma^2}\right) \int_{\mathbb{C}_I} e^{- \frac{p^2}{\gamma^2}}\overline{e_{*}^{2/ \gamma^2}(q \bar{p})} e^{- \frac{\bar{q}^2}{\gamma^2}} e^{\frac{(q- \bar{q})^2}{\gamma^2}} f_I(q) d \lambda_I(q)\\
&=& \left(\frac{2}{\pi \gamma^2}\right) e^{- \frac{p^2}{\gamma^2}} \int_{\mathbb{C}_I} \overline{e_{*}^{2/ \gamma^2}(q \bar{p})} e^{\frac{q^2}{\gamma^2}}  f_I(q) e^{-\frac{2 |q|^2}{\gamma^2}} d \lambda_I(q)\\
&=&\left(\frac{2}{\pi \gamma^2}\right) e^{- \frac{p^2}{\gamma^2}} \int_{\mathbb{C}_I} \overline{e_{*}^{2/ \gamma^2}(q \bar{p})} \mathcal{M}^{\gamma^2}[f_I](q) e^{-\frac{2 |q|^2}{\gamma^2}}d \lambda_I(q).
\end{eqnarray*}
By Theorem \ref{one} we know that $ \mathcal{M}^{\gamma^2}[f](q) \in \mathcal{F}_{Slice}^{\frac{2}{\gamma^2}}(\mathbb{H})$. Then using the reproducing kernel property of the slice hyperholomorphic Fock space, see \eqref{zero}, we get
\begin{eqnarray*}
\left(\frac{2}{\pi \gamma^2}\right) e^{- \frac{p^2}{\gamma^2}} \int_{\mathbb{C}_I} \overline{e_{*}^{2/ \gamma^2}(q \bar{p})} \mathcal{M}^{\gamma^2}[f_I](q) e^{-\frac{2 |q|^2}{\gamma^2}}d \lambda_I(q) &=& e^{- \frac{p^2}{\gamma^2}} \mathcal{M}^{\gamma^2}[f](p).\\
&=& f(p)
\end{eqnarray*}
This concludes the proof.
\end{proof}

\begin{proposition}
	The kernel $K_{\gamma,S}(q,p)$ is slice hyperholomorphic in the variable $q$ and anti-slice hyperholomorphic in variable $p$. 
\end{proposition}
\begin{proof}
The kernel $K_{\gamma,S}(q,p)$ is a pointwise product of an intrinsic slice  hyperholomorphic function $e^{- \frac{q^2}{\gamma^2}}$ and a slice hyperholomorphic function $e_{*}^{\frac{2}{\gamma^2}}(q \bar{p})$. This implies that the kernel $K_{\gamma,S}(q,p)$ is slice hyperholomorphic in the variable $q$, see Theorem \ref{prod}. By similar arguments is anti-slice hyperholomorphic in variable $p$.
\end{proof}

\begin{definition}
\label{ker}
Let $ \gamma >0$. The function
$$ K_{\gamma,S}(q,p)= e^{- \frac{q^2}{\gamma^2}} e_{*}^{\frac{2}{\gamma^2}}(q \bar{p})e^{- \frac{\bar{p}^2}{\gamma^2}},$$
is called the slice hyperholomorphic RBF kernel.
\end{definition}

\begin{proposition}
\label{dkernel}
Let $\gamma >0$. For a fixed $p \in \mathbb{H}$, we set by $K^{p}_{\gamma}$ the function defined as follows
$$ K^{p}_{\gamma,S}(q)= K_{\gamma,S}(q,p).$$
Then it holds that
\begin{itemize}
\item[1)] $K_{\gamma,S}(q,p)= \sum_{n=0}^\infty e_{n}^{\gamma}(q) e_n^\gamma(\bar{p}),$
\item[2)] $ \langle K_{\gamma,S}^q, K_{\gamma,S}^p \rangle_{\mathcal{H}_{\gamma,S}(\mathbb{H})}=K_{\gamma,S}(p,q).$
\end{itemize}
\end{proposition}
\begin{proof}
We start by proving the first point.
\begin{itemize}
\item[1)] Since the functions $q^n$ and $e^{- \frac{q^2}{\gamma^2}}$ are intrinsic slice hyperholomorphic functions they can commute with each other and so by Definition \ref{ker} we have
\begin{eqnarray*}
\sum_{n=0}^\infty e_{n}^\gamma(q) e_n^{\gamma}(\bar{p})&=& \sum_{n=0}^\infty \frac{2^n}{\gamma^{2n} n!} e^{- \frac{q^2}{\gamma^2}} q^n \bar{p}^n e^{- \frac{\bar{p}^2}{\gamma^2}}\\
&=& e^{- \frac{q^2}{\gamma^2}}  e_*^{2/\gamma^2}(q \bar{p})e^{- \frac{\bar{p}^2}{\gamma^2}}\\
&=& K_{\gamma,S}(q,p).
\end{eqnarray*}
\item[2)] For fixed $q$, $p \in \mathbb{H}$ we have that $K_{\gamma}^p(q)$ belongs to the space $\mathcal{H}_{\gamma,S}(\mathbb{H})$. By the reproducing kernel property proved in Theorem \ref{two} we have
$$
\langle K_{\gamma,S}^q, K_{\gamma,S}^p \rangle_{\mathcal{H}_{\gamma,S}(\mathbb{H})}=K^{q}_{\gamma,S}(p)=K_{\gamma,S}(q,p).
$$
\end{itemize}
\end{proof}
For functions in the space $\mathcal{H}_{\gamma,S}(\mathbb{H})$ we have the following estimate.
\begin{proposition}
Let $\gamma>0$ and $I \in \mathbb{S}$. Then for $ f \in \mathcal{H}_{\gamma,S}(\mathbb{H})$ we have
$$ | f(q)| \leq e^{\frac{2}{\gamma^2} y^2} \| f \|_{\mathcal{H}_{\gamma,S}(\mathbb{H})}, \qquad \forall q=x+Iy \in \mathbb{C}_I.$$
\end{proposition}
\begin{proof}
By Theorem \ref{two} we know that
$$ f(q)= \langle f, K_{\gamma,S}^q \rangle_{\mathcal{H}_{\gamma,S}(\mathbb{H})}.$$
By the Cauchy-Schwartz inequality we have
\begin{equation}
\label{f2}
|f(q)| \leq \| K_{\gamma}^q \|_{\mathcal{H}_{\gamma,S}(\mathbb{H})} \| f \|_{\mathcal{H}_{\gamma,S}(\mathbb{H})}.
\end{equation}
By Proposition \ref{dkernel} we get
\begin{eqnarray}
\nonumber
\| K_{\gamma,S}^q \|_{\mathcal{H}_{\gamma,S}(\mathbb{H})}^2 &=& 
\langle K_{\gamma,S}^q, K_{\gamma,S}^q \rangle_{\mathcal{H}_{\gamma,S}(\mathbb{H})}\\
\nonumber
&=&K_{\gamma, S}(q,q)\\
\nonumber
&=& e^{- \frac{q^2}{\gamma^2}} e_{*}^{2/ \gamma^2}(|q|^2) e^{- \frac{\bar{q}^2}{\gamma^2}}\\
\nonumber
&=&e^{- \frac{(q- \bar{q})^2}{\gamma^2}}\\
\label{f3}
&=& e^{\frac{4}{\gamma^2}y^2}.
\end{eqnarray}
By inserting \eqref{f3} in \eqref{f2} we get the thesis.
\end{proof}
A sequential characterization for the slice hyperholomorphic RBF space is also valid. The proof follows similar arguments of \cite[Thm. 3.10]{ADCS}.

\begin{theorem}
A slice hyperholomorphic function, $f(q)= \sum_{n=0}^\infty q^n a_n$, with $ \{a_n\}_{n \geq 0} \subset \mathbb{H}$, belongs to the space $ \mathcal{H}_{\gamma, S}(\mathbb{H})$, if and only if, it holds that
$$ \sum_{k=0}^\infty \frac{k! \gamma^{2k}}{2^k} \left| \sum_{j=0}^{\left[\frac{k}{2}\right]} \frac{a_{k-2j}}{\gamma^{2j} j!}\right|^2< \infty.$$
\end{theorem} 

\begin{proof}
	By Theorem \ref{one} we know that a function $f$ belongs to $ \mathcal{H}_{\gamma,S}(\mathbb{H})$ if and only if there exists a unique function $g \in \mathcal{F}_{Slice}^{\frac{2}{\gamma^2}}(\mathbb{H})$ such that
	$$ f(q)= e^{- \frac{q^2}{\gamma^2}} g(q), \qquad \forall q \in \mathbb{H}.$$
	By \eqref{seq} we know that a function $g(q)= \sum_{k=0}^\infty q^k b_k$, with $ \{b_k\}_{k \geq 0} \subseteq \mathbb{H}$, belongs to $ \mathcal{F}_{Slice}^{\frac{2}{\gamma^2}}(\mathbb{H})$ if we have the following condition
	\begin{equation}
		\label{five}
		\sum_{k=0}^\infty \frac{k! \gamma^{2k}}{2^k} | b_k|^2 < \infty.
	\end{equation}
	By using the Cauchy product of series we have
	\begin{eqnarray}
		\nonumber
		g(q) &=& e^{\frac{q^2}{\gamma^2}} f(q)\\
		\nonumber
		&=& \left( \sum_{n=0}^\infty \frac{q^{2n}}{\gamma^{2n} n!}\right) \left( \sum_{n=0}^\infty q^n a_n\right)\\
		\label{f6}
		&=& \sum_{k=0}^\infty q^k \beta_k
	\end{eqnarray}
	where we set $ \beta_k=\sum_{j=0}^{k} s_j a_{k-j}$ with $s_j=0$ if $j$ is odd, and $s_j=\frac{1}{\gamma^{2m} m!}$ if $j$ is even. Then, for any $ k \geq 0$ we get
	$$
	\beta_k=\sum_{j=0}^k s_j a_{k-j} =\sum_{j=0}^{\left[\frac{k}{2}\right]} \frac{a_{k-2j}}{\gamma^{2j} j!}.
	$$
	By formula \eqref{f6} we have
	$$ g(q)= \sum_{k=0}^\infty q^k b_k= \sum_{k=0}^\infty q^k \beta_k.$$
	By identifying the coefficients of the previous equality we get
	\begin{equation}
		\label{f7}
		b_k= \beta_k =\sum_{j=0}^{\left[\frac{k}{2}\right]} \frac{a_{k-2j}}{\gamma^{2j} j!}, \qquad \forall k \geq 0.
	\end{equation}
	Finally, by plugging \eqref{f7} into \eqref{five} we get
	$$ \sum_{k=0}^\infty \frac{k! \gamma^{2k}}{2^{k}} \left | \sum_{j=0}^{[\frac{k}{2}]} \frac{a_{k-2j}}{\gamma^{2j} j!} \right|^2 < \infty.$$
\end{proof}

\subsection{The slice hyperholomorphic RBF-Segal Bargmann transform }
In this subsection we consider an integral transform that has as a range the slice hyperholomorphic RBF space. This can be connected with the slice hyperholomorphic Bargmann transform. This integral transform was introduced in \cite{DG} and further studied in \cite{DMD}. Now, we recall the main notions of this integral transform. 
The Hilbert space $L^2(\mathbb{R},dx)=L^{2}_{\mathbb{H}}(\mathbb{R})$, which is the domain of the Segal-Bargmann transform, consists of all the square integrable quaternionic-valued functions with respect to 
$$ \langle \varphi, \psi \rangle_{L^2(\mathbb{R},dx)}= \int_{\mathbb{R}} \overline{\psi(x)} \varphi(x) dx.$$
The kernel of the slice hyperholomorphic Segal-Bargmann transform is given by
\begin{equation}
\label{eight}
\mathcal{A}^{\nu}_{SB}(q;x):= \left(\frac{\nu}{\pi}\right)^{\frac{3}{4}} e^{-\frac{\nu}{2}(q^2+x^2)+ \nu \sqrt{2}qx}; \quad (q,x) \in \mathbb{H} \times \mathbb{R},
\end{equation}
see \cite{DG}. This kernel can be seen as the generating function of the real weighted Hermite functions given by
$$ h_n^{\nu}(x):=(-1)^n e^{\frac{\nu^2}{2} x^2} \frac{d}{dx} \left(e^{- \nu x^2}\right),$$
that form an orthogonal basis for the space $L^2(\mathbb{R}, dx)$, with norm given by
$$ \| h_n^{\nu}(x) \|_{L^2(\mathbb{R},dx)}=2^n \nu^n n! \left(\frac{\pi}{\nu}\right)^{\frac{1}{2}}.$$
Namely, we have
\begin{equation}
\label{genker}
\mathcal{A}^{\nu}_{SB}(q,x)=  \sum_{n=0}^\infty  \frac{\sqrt{\nu^n}q^n \psi_n^{\nu}(x)}{\sqrt{n!}},
\end{equation}
where $\psi_n^{\nu}(x)$ are the normalized weighted Hermite functions.

The quaternionic Segal-Bargmann transform $ \mathcal{B}_{\nu}: L^{2}_{\mathbb{H}}(\mathbb{R}) \to \mathcal{F}_{Slice}^{\nu}(\mathbb{H})$ is defined as
$$ \mathcal{B}_{\nu}[\psi](q)= \int_{\mathbb{R}} \mathcal{A}^{\nu}_{SB}(q,x) \psi(x) dx, \qquad \forall \psi \in L^2_{\mathbb{H}}(\mathbb{R}).$$
The action of the quaternionic Segal-Bargmann transform  on the normalized weighted Hermite functions is given by
\begin{equation}
\label{action}
 \mathcal{B}_{\nu}[\psi_n^{\nu}](q)= \left(\frac{\nu}{\pi}\right)^{\frac{1}{2}} \frac{\nu^{\frac{n}{2}} q^n}{\sqrt{n!}},
\end{equation}
see \cite[Lemma 4.4]{DG}. We can define the kernel of the slice hyperholomorphic RBF-Segal Bargmann transform as
$$ \mathcal{A}_{RBF}^{\gamma}(q,x)= \sum_{n=0}^\infty e_n^{\gamma}(q) \psi_{n}^{\frac{2}{\gamma^2}}(x), \qquad (q,x) \in \mathbb{H} \times \mathbb{R}.$$ 
\begin{definition}
The slice hyperholomorphic RBF-Segal Bargmann transform is defined as
$$ \mathcal{B}_{RBF}^{\gamma}[\psi](q)=\int_{\mathbb{R}} \mathcal{A}_{RBF}^{\gamma}(q,x) \psi(x) dx,$$
for any $\psi \in L^{2}_{\mathbb{H}}(\mathbb{R})$ and $q\in \mathbb{H}$.
\end{definition}

Now, we show the following match between the kernels of the slice hyperholomorphic Segal-Bargmann and the slice hyperholomorphic RBF-Segal-Bargmann transforms. 

\begin{proposition}
\label{twobis}
Let $\gamma>0$. Then we have
$$ \mathcal{A}_{RBF}^{\gamma}(q,x)=e^{-\frac{q^2}{\gamma^2}}  \mathcal{A}_{SB}^{\frac{2}{\gamma^2}}(q,x), \quad \forall (q,x) \in \mathbb{H} \times \mathbb{R}.$$
\end{proposition}
\begin{proof}
By the kernel defined in \eqref{genker} and the definition of the basis $e_{n}^\gamma(q)$ (see \eqref{basis}) we get
$$
\mathcal{A}_{RBF}^{\gamma}(q,x)=e^{- \frac{q^2}{\gamma^2}} \sum_{n=0}^\infty \sqrt{\frac{2^n}{n! \gamma^{2n}} } q^n \psi_{n}^{\frac{2}{\gamma^2}}(x)=e^{- \frac{q^2}{\gamma^2}} \mathcal{A}_{SB}^{\frac{2}{\gamma^2}}(q,x).
$$
\end{proof}
A relation  between the slice hyperholomorphic Segal-Bargmann transform and the RBF one is possible. 
\begin{corollary}
\label{four}
Let $\gamma>0$ and $\psi \in L^2_{\mathbb{H}}(\mathbb{R})$. Then we have
$$ \mathcal{B}^{\gamma}_{RBF}[\psi](q)= \mathcal{M}^{- \gamma^2} \mathcal{B}_{\frac{2}{\gamma^2}} [\psi](q), \quad \forall q \in \mathbb{H}.$$
\end{corollary}
\begin{proof}
The result follows by the definition of the slice hyperholomorphic RBF-Segal-Bargmann transform and Proposition \ref{twobis}. 
\end{proof}

\begin{proposition}
The explicit expression of the slice hyperholomorphic RBF-kernel is given by
$$ \mathcal{A}_{RBF}^{\gamma}(q,x)= \left(\frac{2}{\pi \gamma^2}\right)^{\frac{3}{4}} e^{- \frac{(x- \sqrt{2}q)^2}{\gamma^2}} \qquad \forall (q,x) \in \mathbb{H} \times \mathbb{R}.$$ 
\end{proposition}
\begin{proof}
By Proposition \ref{twobis} we have
\begin{eqnarray*}
\mathcal{A}_{RBF}^{\gamma}(q,x)&=& e^{- \frac{q^2}{\gamma^2}} \mathcal{A}_{SB}^{\frac{2}{\gamma^2}}(q,x)\\
&=& \left( \frac{2}{ \pi \gamma^2}\right)^{\frac{3}{4}} e^{- \frac{2q^2}{\gamma^2}- \frac{x^2}{\gamma^2}+ \frac{2 \sqrt{2} qx}{\gamma^2}}\\
&=&  \left( \frac{2}{ \pi \gamma^2}\right)^{\frac{3}{4}} e^{-\frac{(x- \sqrt{2}q)^2}{\gamma^2}}.
\end{eqnarray*}
\end{proof}
\begin{proposition}
Let $\gamma > 0$ and $n  \in \mathbb{N}$. Then we have
$$ \mathcal{B}_{RBF}^{\gamma}[\psi_n^{\frac{2}{\gamma^2}}](q)= \sqrt{\frac{2}{\gamma^2 \pi}}e_{n}^{\gamma}(q).$$
Moreover, we have also
$$ \|  \mathcal{B}_{RBF}^{\gamma}[\psi_n^{\frac{2}{\gamma^2}}] \|_{\mathcal{H}_{\gamma,S}(\mathbb{H})}= \| \psi_{n}^{\frac{2}{\gamma^2}} \|_{L^{2}_{\mathbb{H}}(\mathbb{R})}=1.$$
\end{proposition}
\begin{proof}
By Corollary \ref{four} and by formula \eqref{action} with $ \nu = \frac{2 }{\gamma^2}$ we get
\begin{eqnarray*}
\mathcal{B}^{\gamma}_{RBF}[\psi_n^{\frac{2}{\gamma^2}}](q)&=& \mathcal{M}^{- \gamma^2}\mathcal{B}_{\frac{2}{\gamma^2}}[\psi_n^{\frac{2}{\gamma^2}}](q)\\
&=& e^{- \frac{q^2}{\gamma^2}} \mathcal{B}_{\frac{2}{\gamma^2}}[\psi_n^{\frac{2}{\gamma^2}}](q)\\
&=& \left( \frac{2}{\gamma^2 \pi}\right)^{\frac{1}{2}} \sqrt{\frac{2^n}{\gamma^{2n} n!}} q^n e^{- \frac{q^2}{\gamma^2}}\\
&=& \sqrt{\frac{2}{\gamma^2 \pi}} e_{n}^{\gamma}(q).
\end{eqnarray*}
\end{proof}

\begin{theorem}
For $\gamma > 0$ the slice hyperholomorphic RBF Bargmann transform is an isometric isomorphism mapping $L^2_{\mathbb{H}}(\mathbb{R})$ onto $ \mathcal{H}_{\gamma,S}(\mathbb{H})$. 
\end{theorem}
\begin{proof}
By Corollary \ref{four}, Theorem \ref{zero} and \cite[Thm. 4.6]{DG} we get
\begin{eqnarray*}
\| \mathcal{B}_{RBF}^{\gamma}[\psi] \|_{\mathcal{H}_{\gamma,S}(\mathbb{H})}&=&\|\mathcal{M}^{- \gamma^2}\mathcal{B}_{\frac{2}{\gamma^2}}[\psi] \|_{\mathcal{H}_{\gamma,S}(\mathbb{H})}\\
&=& \| \mathcal{B}_{\frac{2}{\gamma^2}}[\psi]\|_{\mathcal{F}_{Slice}^{\frac{2}{\gamma^2}}}\\
&=& \| \psi \|_{L^2_{\mathbb{H}}(\mathbb{R})},
\end{eqnarray*}
for any $\psi \in L^{2}_{\mathbb{H}}(\mathbb{R})$.
\end{proof}

\section{Gaussian RBF-kernels via Fock spaces: Several complex variables case}

In this section we extend some results obtained in \cite{ADCS} to the case of several complex variables. We present some of the results without proofs since they follow similar arguments of the one complex variable. To this end, let us first introduce some standard notations that are used in the case of several complex variables with dimension $d\geq 1$. On the space $\mathbb{C}^d$ we can consider the inner product and the associated metric given by 
\begin{equation}
	\langle z,w \rangle_{2,d}=\displaystyle \sum_{\ell=1}^{d}z_\ell \overline{w_\ell}, \quad\displaystyle d(z,w)=||z-w||_{2,d}=\sqrt{\sum_{\ell=1}^{d}|z_\ell-w_\ell|^2},
\end{equation}
for every $z=(z_1,...,z_d), w=(w_1,...,w_d)\in\mathbb{C}^d$. We use also the following notations for the conjugate and product
\begin{equation}\label{Notaation}
	\displaystyle \overline{z}=(\overline{z}_1,\cdots,\overline{z}_d), \quad z\cdot w= \sum_{\ell=1}^{d} z_\ell w_\ell, \text{ and }  \quad z\cdot z=z^2=\sum_{\ell=1}^{d} z^2_\ell.
\end{equation}
Given $n=(n_1,\cdots,n_d)\in \mathbb{N}^d$ and $z=(z_1,...,z_d)\in \mathbb{C}^d$ we use the multi-index notation: 
$$z^n=z_{1}^{n_1}z_{2}^{n_2}\cdots z_{d}^{n_d}, \quad n!=(n_1!)(n_2!)\cdots (n_d!), \quad |n|=\sum_{\ell=1}^{d}n_\ell.$$

Let $\alpha>0$. We briefly recall the Fock space on $\mathbb{C}^d$ which is a subspace of entire functions $Hol(\mathbb{C}^d)$ defined by 

$$\mathcal{F}_\alpha(\mathbb{C}^d)=\left\lbrace f\in Hol(\mathbb{C}^d), \quad \left(\frac{\alpha}{\pi}\right)^d\int_{\mathbb{C}^d}|f(z)|^2e^{-{\alpha}|z|^2}dA(z)<\infty \right\rbrace, $$
where $dA(z)$ denotes the Lebesgue measure on $\mathbb{C}^d$ with $z=(z_1,\cdots,z_d)\in \mathbb{C}^d$. The reproducing kernel of the space $\mathcal{F}_\alpha(\mathbb{C}^d)$ is given by 
\begin{equation}
	F_\alpha(z,w)=e^{ \alpha z\cdot \overline{w}}, \quad \forall z=(z_1,...,z_d), w=(w_1,...,w_d)\in \mathbb{C}^d.
\end{equation}

The authors of \cite{SDC} (see also the book \cite{SC2008}) introduced the complex Gaussian RBF kernel in several complex variables as follows: 

\begin{definition}
	Let $d\geq 1$  and $\gamma >0$. The $d$-dimensional complex valued Gaussian RBF kernel on $\mathbb{C}^d$ is denoted by $K_{\gamma,d}: \mathbb{C}^d\times \mathbb{C}^d\longrightarrow \mathbb{C}$ and can be expressed as follows
	$$K_{\gamma,d}(z,w)=\exp\left(-\dfrac{(z-\overline{w})^2}{\gamma^2}\right), $$
	for every $ z=(z_1,...,z_d),w=(w_1,...w_d)\in\mathbb{C}^d.$
\end{definition}

\begin{remark}
	If we restrict the complex variables $z=(z_1,...,z_d),w=(w_1,...w_d)\in\mathbb{C}^d$ to the real variables $x=(x_1,...,x_d),y=(y_1,...y_d)\in\mathbb{R}^d$ we re-obtain the classical real valued Gaussian RBF kernel
	
$$
		K_{\gamma,d}(x,y)=\exp\left(-\dfrac{||x-y||^{2}_{2,d}}{\gamma^2}\right).
$$
\end{remark}

We make the following observation:

\begin{proposition}\label{Prop3}
	Let $d\geq 1$ and $\gamma>0$. The Gaussian RBF kernel of dimension $d$ can be expressed as a product of Gaussian RBF kernels of dimension one. More precisely, we have 
	\begin{equation}
		\label{new1}
		K_{\gamma,d}(z,w)= \prod_{\ell=1}^{d} K_\gamma(z_\ell,w_\ell),
	\end{equation}
	for every $z=(z_1,...,z_d), w=(w_1,...,w_d)\in\mathbb{C}^d$.
\end{proposition}

\begin{proof}
	Let $d\geq 1$ and $\gamma>0$. For $z=(z_1,...,z_d)$ and $w=(w_1,...,w_d)$ in $\mathbb{C}^d$ we have
	\[ \begin{split}
		\displaystyle \prod_{\ell=1}^{d} K_\gamma(z_\ell,w_\ell) & = \prod_{\ell=1}^{d}\exp\left(-\dfrac{(z_\ell-\overline{w}_\ell)^2}{\gamma^2}\right) \\
		& =\exp\left(-\frac{1}{\gamma^2}\sum_{\ell=1}^{d}(z_\ell-\overline{w}_\ell)^2\right) \\
		& =\exp\left(-\frac{(z-\overline{w})^2}{\gamma^2}\right)\\
		&=K_{\gamma,d}(z,w).\\
	\end{split}
	\]     
	
\end{proof}

\begin{remark}
The function $K_{\gamma,d}(z,w)$ (see \eqref{new1}) is a kernel as a product of one dimensional RBF kernels, see \cite[Lemma 4.6]{SC2008}. Moreover, using the tensor product symbol $\otimes$, we observe the following 
	$$K_{\gamma,d}(z,w)=\exp\left(-\dfrac{(z-\overline{w})^2}{\gamma^2}\right)=K_{\gamma}^{\otimes d}(z,w), $$
	for every $ z=(z_1,...,z_d),w=(w_1,...w_d)\in\mathbb{C}^d.$
\end{remark}

\begin{proposition}
	The complex Gaussian RBF kernel $K_\gamma$ can be expressed in terms of the reproducing kernel of the Fock space $\mathcal{F}_{\frac{2}{\gamma^2}}(\mathbb{C}^d)$. More precisely, we have   
$$
		K_{\gamma,d}(z,w)=\exp\left(-\frac{1}{\gamma^2}(z^2+\overline{w}^2)\right)F_{\frac{2}{\gamma^2}}(z,w),
$$
	for every $z=(z_1,\cdots,z_d)$ and $w=(w_1,\cdots,w_d)$ in $\mathbb{C}^d$.
\end{proposition}
\begin{proof}
	We know by Proposition \ref{Prop3} that the Gaussian RBF kernel of dimension $d$ can be expressed as a product of Gaussian RBF kernels of dimension $1$. Then, by \cite[Proposition 3.1]{ADCS} we get 
	\begin{align*}
		K_{\gamma,d}(z,w)&=\prod_{\ell=1}^{d} K_\gamma(z_\ell,w_\ell)\\
		&=\prod_{\ell=1}^{d}\exp\left(-\frac{(z_\ell^2+\overline{w_\ell}^2)}{\gamma^2}\right)F_{\frac{2}{\gamma^2}}(z_\ell,w_\ell)\\
		&=\exp\left(-\frac{\sum_{\ell=1}^{d}(z_\ell^2+\overline{w_\ell}^2)}{\gamma^2}\right)F_{\frac{2}{\gamma^2}}(z,w)\\
		&=\exp\left(-\frac{1}{\gamma^2}(z^2+\overline{w}^2)\right)F_{\frac{2}{\gamma^2}}(z,w).\\
	\end{align*}
\end{proof}

The complex reproducing kernel Hilbert space associated to the Gaussian RBF kernel in the case of several complex variables can be introduced as follows (see \cite{SDC, SC2008}):

\begin{definition}
	Let $\gamma>0$, an entire function $f:\mathbb{C}^d\longrightarrow \mathbb{C}$ belongs to the RBF space, denoted by $\mathcal{H}_{\gamma}(\mathbb{C}^d)$ (or simply $\mathcal{H}_{\gamma,d})$ if we have
	$$
	\displaystyle ||f||_{\mathcal{H}_{\gamma,d}}^2:=\left(\dfrac{2}{\pi\gamma^2}\right)^d\int_{\mathbb{C}^d}|f(z)|^2\exp\left(\frac{(z-\overline{z})^2}{\gamma^2}\right) dA(z)<\infty,$$
	where $dA(z)=\prod_{\ell=1}^{d} dx_\ell dy_\ell$ is the Lebesgue measure with respect to the variable $z=x+iy\in \mathbb{C}^d$.
\end{definition}
We can apply similar arguments used in \cite{ADCS} in the case of one complex variable to prove the following result for RBF spaces on $\mathbb{C}^d$:
\begin{theorem}\label{RKP-RBF}
	Let $\gamma>0$, the RBF Hilbert space $\mathcal{H}_{\gamma,d}$ is a reproducing kernel Hilbert space whose reproducing kernel is given by the RBF kernel $K_{\gamma,d}(z,w)$. Moreover, we have the reproducing kernel property which is given by the following integral representation
	$$\label{RPRBF}
\displaystyle f(w)=\left(\dfrac{2}{\pi\gamma^2}\right)^d\int_{\mathbb{C}^d} f(z)\overline{K_{\gamma,d}(z,w)}\exp\left(\frac{(z-\overline{z})^2}{\gamma^2}\right)dA(z),
	$$
	for any $ f\in\mathcal{H}_{\gamma,d}, w=(w_1,\cdots,w_d)\in\mathbb{C}^d$. 
\end{theorem}
\begin{remark}
	In other terms, we can reproduce any function $f$ belonging to the complex $d$-dimensional RBF space $\mathcal{H}_{\gamma,d}$ using the following reproducing kernel property 
	$$f(w)=\langle f, K_{\gamma,d}^{w} \rangle_{\mathcal{H}_{\gamma,d}},\qquad f\in\mathcal{H}_{\gamma,d}, \quad w=(w_1,\cdots,w_d)\in\mathbb{C}^d.$$
\end{remark}

It is important to note that an orthonormal basis of the Gaussian RBF space $\mathcal{H}_{\gamma,d}$ is given by the family of functions defined for every $n=(n_1,\cdots,n_d)\in \mathbb{N}^d$ by the following expression (see \cite{SDC, SC2008}):
\begin{equation}
	e_{n}^{\gamma}(z)= \sqrt{\frac{2^{|n|}}{\gamma^{2|n|} n!}} z^n e^{-\frac{z^2}{\gamma^2}}, \qquad \forall z=(z_1,\cdots,z_d)\in \mathbb{C}^d.
\end{equation}
We can prove the following expression of the Gaussian RBF kernel:
\begin{proposition}
	For every $z=(z_1,\cdots,z_d), w=(w_1,\cdots,w_d)\in \mathbb{C}^d$ we have
$$
		\displaystyle \sum_{n=(n_1,\cdots,n_d)\in \mathbb{N}^d} e_{n}^{\gamma}(z)\overline{e_{n}^{\gamma}(w)}=K_{\gamma,d}(z,w).
$$
\end{proposition}
\begin{proof}
	We apply the one dimensional case (see \cite[Proposition 3.5]{ADCS}) and get 
	\begin{align*}
		\displaystyle \sum_{n=(n_1,\cdots,n_d)\in \mathbb{N}^d} e_{n}^{\gamma}(z)\overline{e_{n}^{\gamma}(w)}&=\prod_{\ell=1}^{d}\left( \sum_{n_\ell=0}^{\infty} e_{n_\ell}^{\gamma}(z_\ell)\overline{e_{n_\ell}^{\gamma}(w_\ell)}\right)\\
		&=\prod_{\ell=1}^{d}K_\gamma(z_\ell,w_\ell)\\
		&=K_{\gamma,d}(z,w).\\
	\end{align*}
\end{proof}
Inspired from Theorem 3.2 of \cite{ADCS} we note that the Gaussian RBF space of several complex variables $\mathcal{H}_{\gamma,d}$ is isomoprhic to the Fock space on $\mathbb{C}^d$ with parameter $\frac{2}{\gamma^2}$. More precisely, we have the following result 
\begin{theorem}[RBF-Fock isomorphism]\label{FRBFchara}
	Let $\gamma>0$, an entire function $f:\mathbb{C}^d\longrightarrow \mathbb{C}$ belongs to the RBF space $\mathcal{H}_{\gamma,d}$ if and only if there exists a unique function $g$ in the Fock space $\mathcal{F}_{\frac{2}{\gamma^2}}(\mathbb{C}^d)$such that $$f(z)=\exp(-\frac{z^2}{\gamma^2})g(z), \quad \text{ for any } z=(z_1,\cdots,z_d)\in\mathbb{C}^d.$$ Moreover, there exists an isometric isomorphism between the RBF and Fock spaces given by the multiplication operator $\mathcal{M}_{RBF}^{\gamma^2}:\mathcal{H}_{\gamma,d}\longrightarrow \mathcal{F}_{\frac{2}{\gamma^2}}(\mathbb{C}^d)$ defined by
$$
		\mathcal{M}^{\gamma^2}_{RBF}[f](z):=\mathcal{M}_{\exp{(\frac{z^2}{\gamma^2}})}[f](z)=\exp(\frac{z^2}{\gamma^2})f(z), \quad \text{ for any } f\in \mathcal{H}_{\gamma,d}, \quad z\in \mathbb{C}^d.
$$
\end{theorem}
Finally, the counterpart of the RBF-Segal Bargmann transform for several complex variables can be introduced for any $\varphi \in L^2(\mathbb{R}^d)$ and $z=(z_1,\cdots,z_d)\in\mathbb{C}^d$ as follows:

	$$\displaystyle \mathfrak{B}_{\gamma,d}[\varphi](z_1,\cdots,z_d)=\int_{\mathbb{R}^d} \mathfrak{A}_{RBF}^{\gamma,d}((z_1,\cdots,z_d), (x_1,\cdots,x_d))\varphi(x_1,\cdots,x_d)dx_1\cdots dx_d,$$
	where the RBF-Bargmann kernel here is obtained by taking the generating function associated to the normalized Hermite functions $\psi_n(x_1,\cdots,x_d)=\psi_{n_1}(x_1)\cdots \psi_{n_d}(x_d)$ and the orthonormal basis $e_n^\gamma(z_1,\cdots,z_d)$ with $n=(n_1,\cdots,n_d)$ so that:
	
	$$\displaystyle \mathfrak{A}_{RBF}^{\gamma,d}((z_1,\cdots,z_d), (x_1,\cdots,x_d))=\sum_{n\in \mathbb{N}^d}\psi_n(x_1,\cdots,z_d)e_n(z_1,\cdots,z_d).$$
By following a similar approach used in Proposition \ref{twobis} we have
	
	$$\sum_{n\in \mathbb{N}^d}\psi_n(x_1,\cdots,z_d)e_n(z_1,\cdots,z_d)=\left(\dfrac{2}{\pi\gamma^2 }\right)^{\frac{d}{4}} \exp\left(-\frac{\sum_{\ell=1}^{d}(\sqrt{2}z_\ell-x_\ell)^2}{\gamma^2}\right). $$
	As a consequence it is possible to express the RBF-Segal-Bargmann transform as follows
	
	\begin{equation}\label{SBT}
		\displaystyle \mathfrak{B}_{\gamma,d}[\varphi](z_1,\cdots,z_d)=\left(\dfrac{2}{\pi\gamma^2 }\right)^{\frac{d}{4}}\int_{\mathbb{R}^d} \exp\left(-\frac{\sum_{\ell=1}^{d}(\sqrt{2}z_\ell-x_\ell)^2}{\gamma^2}\right)\varphi(x_1,\cdots,x_d)dx_1\cdots dx_d.
\end{equation} 
We observe that using the standard notations in \eqref{Notaation} we can rewrite the previous expression in the following way
$$
\displaystyle \mathfrak{B}_{\gamma,d}[\varphi](z)=\left(\dfrac{2}{\pi\gamma^2 }\right)^{\frac{d}{4}}\int_{\mathbb{R}^d} \exp\left(-\frac{(\sqrt{2}z-x)^2}{\gamma^2}\right)\varphi(x)dx.
$$

The following result holds true 
\begin{theorem}
	The RBF Segal-Bargmann transform defined by the expression \eqref{SBT}
	is an isometric isomorphism mapping the standard Hilbert space $L^2(\mathbb{R}^d)$ onto the RBF space $\mathcal{H}_{\gamma,d}$.
\end{theorem}

\section{Concluding remarks}
In this paper we have investigated Gaussian RBF kernels and their connections to the theory of Fock spaces in two settings: quaternionic slice hyperholomorphic and several complex  variables. In a forthcoming paper we aim to study a generalization of this topic in the case of quaternionic monogenic functions using the Fueter mapping theorem, which allows to construct monogenic functions starting from slice hyperholomorphic ones. The various extensions of RBF kernel are summarized in the following scheme:
{\small
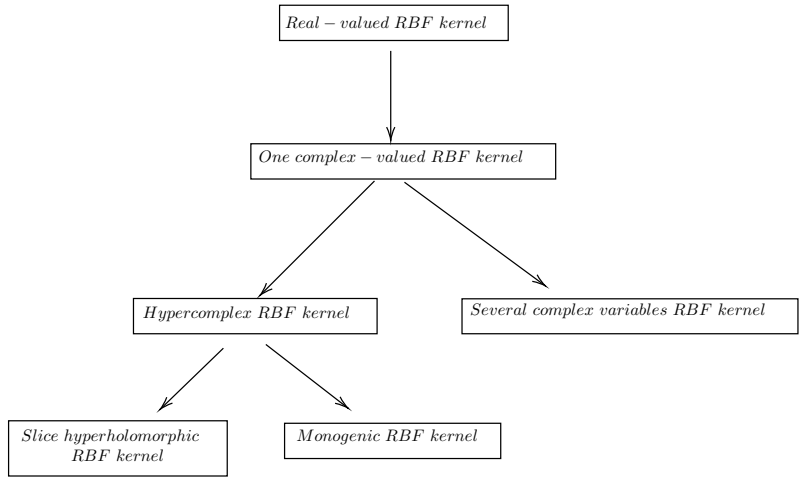
	\begin{figure}[H]
		\centering
		\resizebox{0.65\textwidth}{!}{%
			\tikzset{every picture/.style={line width=0.75pt}} 

\begin{tikzpicture}[x=0.75pt,y=0.75pt,yscale=-1,xscale=1]
	
	\draw    (327,160.5) -- (236.5,251) ;
	\draw [shift={(235.08,252.42)}, rotate = 315] [color={rgb, 255:red, 0; green, 0; blue, 0 }  ][line width=0.75]    (10.93,-3.29) .. controls (6.95,-1.4) and (3.31,-0.3) .. (0,0) .. controls (3.31,0.3) and (6.95,1.4) .. (10.93,3.29)   ;
	\draw    (206,294.5) -- (156.52,342.36) ;
	\draw [shift={(155.08,343.75)}, rotate = 315.95] [color={rgb, 255:red, 0; green, 0; blue, 0 }  ][line width=0.75]    (10.93,-3.29) .. controls (6.95,-1.4) and (3.31,-0.3) .. (0,0) .. controls (3.31,0.3) and (6.95,1.4) .. (10.93,3.29)   ;
	\draw    (240,291.5) -- (304.5,341.52) ;
	\draw [shift={(306.08,342.75)}, rotate = 217.79] [color={rgb, 255:red, 0; green, 0; blue, 0 }  ][line width=0.75]    (10.93,-3.29) .. controls (6.95,-1.4) and (3.31,-0.3) .. (0,0) .. controls (3.31,0.3) and (6.95,1.4) .. (10.93,3.29)   ;
	\draw    (351,161.5) -- (462.47,243.56) ;
	\draw [shift={(464.08,244.75)}, rotate = 216.36] [color={rgb, 255:red, 0; green, 0; blue, 0 }  ][line width=0.75]    (10.93,-3.29) .. controls (6.95,-1.4) and (3.31,-0.3) .. (0,0) .. controls (3.31,0.3) and (6.95,1.4) .. (10.93,3.29)   ;
	\draw    (340,56.25) -- (340,126.25) ;
	\draw [shift={(340,128.25)}, rotate = 270] [color={rgb, 255:red, 0; green, 0; blue, 0 }  ][line width=0.75]    (10.93,-3.29) .. controls (6.95,-1.4) and (3.31,-0.3) .. (0,0) .. controls (3.31,0.3) and (6.95,1.4) .. (10.93,3.29)   ;
	\draw   (251,19.67) -- (434.08,19.67) -- (434.08,47.92) -- (251,47.92) -- cycle ;
	\draw   (228,130.67) -- (472.08,130.67) -- (472.08,158.92) -- (228,158.92) -- cycle ;
	\draw   (397.08,254.92) -- (666.08,254.92) -- (666.08,282.92) -- (397.08,282.92) -- cycle ;
	\draw   (134,254.67) -- (329.08,254.67) -- (329.08,282.92) -- (134,282.92) -- cycle ;
	\draw   (255,354.67) -- (428.08,354.67) -- (428.08,382.92) -- (255,382.92) -- cycle ;
	\draw   (34,352.67) -- (209.08,352.67) -- (209.08,396.92) -- (34,396.92) -- cycle ;
	
	\draw (233,134.9) node [anchor=north west][inner sep=0.75pt]    {$One\ complex-valued\ RBF\ kernel$};
	\draw (140,259.9) node [anchor=north west][inner sep=0.75pt]    {$Hypercomplex\ RBF\ kernel$};
	\draw (402,256.9) node [anchor=north west][inner sep=0.75pt]    {$Several\ complex\ variables\ RBF\ kernel$};
	\draw (36,356.07) node [anchor=north west][inner sep=0.75pt]    {$ \begin{array}{l}
			Slice\ hyperholomorphic\\
			\ \ \ \ \ \ \ \ \ RBF\ kernel
		\end{array}$};
	\draw (263,357.9) node [anchor=north west][inner sep=0.75pt]    {$Monogenic\ RBF\ kernel$};
	\draw (254,27.4) node [anchor=north west][inner sep=0.75pt]    {$Real-valued\ RBF\ kernel$};

\end{tikzpicture}
		}
	\end{figure}
}

\textbf{Acknowledgements}
The authors are grateful to the anonymous referee whose deep and extensive comments greatly contributed to improve this paper.


\begin{thebibliography}{AA}
\bibitem{ADCS} D.~Alpay, F.~Colombo, K.~Diki, I.~Sabadini, \emph{An approach to the Gaussian RBF kernels via Fock spaces}. J. Math. Phys. \textbf{63} (2022), no. 11, Paper No. 113506, 19 pp.
\bibitem{ACSbook} D.~Alpay , F.~Colombo, I.~Sabadini, \emph{Slice hyperholomorphic Schur analysis. Operator Theory}: Advances and Applications, 256. Birkhäuser/Springer, Cham, 2016. xii+362 pp.
\bibitem{6COFBook}
D.~Alpay, F.~Colombo,  I.~Sabadini,
{\em Quaternionic de Branges spaces and characteristic operator function},
SpringerBriefs in Mathematics, Springer, Cham, 2020/21
\bibitem{ACSS} D.~Alpay, F.~Colombo, I.~Sabadini, G.~Salomon, \emph{ The Fock space in the slice hyperholomorphic Setting}. In: Hypercomplex Analysis: New perspectives and applications. Trends Math. pp. 43–59(2014).
\bibitem{red}  F.~Brackx , R.~Delanghe, F.~Sommen, {\em Clifford Analysis}, Pitman Res. Notes in Math., 76, 1982.



\bibitem{frac5}
F. Colombo, D. Deniz-Gonzales, S. Pinton,
{\em Non commutative fractional Fourier law in bounded and unbounded domains},
Complex Analyis and Operator Theory 2021.

\bibitem{6CG}
F. Colombo, J. Gantner,
{\em Quaternionic closed operators, fractional powers and fractional diffusion processes},
Operator Theory: Advances and Applications, 274. Birkh\"auser/Springer, Cham, 2019. viii+322 pp.
\bibitem{CGK} F.~Colombo, J.~Gantner, P.~Kimsey, \emph{Spectral theory on the S-spectrum for quaternionic operators}. In: Operator Theory: Advances and Applications, 270. Birkhauser/Springer, Cham, p. ix+356 (2018).


\bibitem{frac3}
F. Colombo,  M. Peloso, S. Pinton,
{\em The structure of the fractional powers of the noncommutative Fourier law},
Math. Methods Appl. Sci., {\bf 42} (2019),  6259--6276
\bibitem{CSS1} F.~Colombo , I.~Sabadini, D.C.~Struppa, \emph{Noncommutative functional calculus, Progress in Mathematics}, vol. 289, Birkhäuser/Springer Basel AG, Basel, (2011).
\bibitem{CSS2}F.~Colombo , I.~Sabadini I., D.C.~Struppa, \emph{Entire Slice Regular Functions}, SpringerBriefs in Mathematics, Springer, Cham, (2016).
\bibitem{CSS3} F.~Colombo, I.~Sabadini, D.~C.~Struppa, \emph{Michele Sce's Works in Hypercomplex Analysis. A Translation with Commentaries}, Birkhäuser/Springer Basel AG, Basel, 2020.
\bibitem{DMD} A.~De Martino, K.~Diki, \emph{On the Quaternionic Short-Time Fourier and Segal–Bargmann Transforms},  Mediterr. J. Math. \textbf{18}, 110 (2021).
\bibitem{DMD2} A.~De Martino , K.~Diki, \emph{On the polyanalytic short-time Fourier tranform in the quaternionic setting}, Commun. Pure Appl. Anal, 21, \textbf{11}, (2022), 3629–3665.
\bibitem{GS2007}  G.~Gentili , C.~Stoppato , D.C.~Struppa , \emph{Regular functions of a quaternionic varaible}, Springer Monographs, Berlin (2013).
\bibitem{F} R.~Fueter, \emph{ Die {F}unktionentheorie der {D}ifferentialgleichungen $\Delta u=0$ und {$\Delta\Delta u=0$} mit vier reellen {V}ariablen}, Comment. Math. Helv. 7(1):307--330, 1934.
\bibitem{DG} K.~Diki, A.~Ghanmi, \emph{A quaternionic analogue for the Segal-Bargamann transfrom}, Complex. Anal. Oper. Theory \textbf{11}, 457-473 (2017).
\bibitem{Gaz}  J.P.~Gazeau, \emph{Coherent states in quantum physics}, WILEY-VCH Verlag GmbH and Co. KGaA, Weinheim (2009).
\bibitem{SHS} B. Schölkopf, R. Herbrich, AJ. Smola, \emph{A generalized representer theorem}, In Computational Learning Theory: 14th Annual Conference on Computational Learning Theory, COLT 2001 and 5th European Conference on Computational Learning Theory, EuroCOLT 2001 Amsterdam, The Netherlands, July 16–19, 2001 Proceedings 14 2001 (pp. 416-426). Springer Berlin Heidelberg.
\bibitem{SDC}I.~Steinwart, H.~Don, S.~Clint, \emph{ An explicit description of the reproducing kernel Hilbert spaces of Gaussian RBF kernels}, IEEE Transactions on Information Theory 52.10 (2006): 4635-4643.
\bibitem{SC2008}I.~Steinwart, A. Christmann, Support Vector Machines (Springer-Verlag New York, 2008), ISBN: 978-0-387-77242-4

\bibitem{vert2004primer}
J.P. Vert, K. Tsuda, and B. Sch{\"o}lkopf.
\newblock A primer on kernel methods.
\newblock { Kernel methods in computational biology}, 47:35--70, 2004.


	\end{thebibliography}
\end{document}